\documentclass{article}
\usepackage{amsmath,amsfonts,epsfig,xcolor}
\usepackage{mathptmx}
\usepackage{amsmath,amsfonts,latexsym,amssymb}
\usepackage{mathrsfs}
\usepackage{verbatim}
\usepackage[utf8]{inputenc}
\usepackage[T1]{fontenc}
\usepackage{ae,aecompl}
\usepackage{braket}

\newcommand{\bb}{\mathbb}

\newcommand{\cx}{{\bb C}}

\renewcommand{\bold}[1]{\medskip \noindent {\bf \boldmath #1
                        }\nopagebreak[4]}

\newcommand{\del}{\partial}

\newcommand{\chat}{\widehat{\cx}}

\newcommand{\area}{\operatorname{area}}

\newcommand{\tr}{\operatorname{tr}}
\newcommand{\vol}{\operatorname{vol}}

\newtheorem{theorem}{Theorem}[section]
\newtheorem{prop}[theorem]{Proposition}
\newtheorem{lemma}[theorem]{Lemma}

\newcommand{\cB}{{\cal B}}

\newcommand{\cM}{{\cal M}}

\newtheorem{corollary}[theorem]{Corollary}

\renewcommand{\hbar}{\bar{{\mathbb H}}^3}

\newcommand{\CC}{\mathbb C}

\newcommand{\Hs}{{{\mathbb H}^3}}

\newcommand{\psl}{{\rm PSL}(2,\CC)}

\def\eproof{$\Box$ \medskip}
\renewcommand{\area}{{\operatorname{{\bf area}}}}

\newcommand{\Ep}{\operatorname{Ep}}

\newcommand{\Id}{{\operatorname{Id}}}

\begin{document}

\title{\bf A bound on the $L^2$-norm of a projective structure by the length of the  bending lamination}
\author{Martin
   Bridgeman\thanks{M. Bridgeman's research  supported by  NSF grant DMS-2005498.} \ and
  Kenneth Bromberg\thanks{K. Bromberg's research supported by NSF grant
      DMS-1906095.}}

\date{\today}

\maketitle

\begin{abstract}
One can associate to a complex projective structure on a surface holomorphic quadratic differential $\Phi$ via the  Schwarzian derivative and a bending lamination $\lambda$ via the Thurston parameterization. In this note we obtain upper bounds on the $L^2$-norm of $\Phi$ in terms of the length of $\lambda$. The proof uses the theory of $W$-volume introduced by Krasnov-Schlenker.
\end{abstract}


\section{Introduction}

Complex projective structures on a surface $S$ have two parameterizations. Classically they are given via the Schwarzian derivative by a pair $(X,\Phi)$ where $X$ is a conformal structure on $S$ and $\Phi$ is a holomorphic quadratic differential on $X$. Thurston gave an alternative parameterization by pairs $(Y,\lambda)$ where $Y$ is a hyperbolic structure and $\lambda$ is a measured geodesic lamination on $Y$. Both of these parameterizations give a way of measuring the ``size'' of the projective structure. The space of holomorphic quadratic differentials on $X$ is similar to a function space and we can take the $L^p$-norms of $\Phi$ with respect to the hyperbolic metric on $X$ while a measured lamination on a hyperbolic surface has a length $L(\lambda)$. Naturally one would like to obtain comparisons between the norm of $\Phi$ and the length $L(\lambda)$. Note that the subspace of projective structures where $\|\Phi\|_p = 0$ is equal to the subspace where $L(\lambda)=0$ and  this happens exactly when the conformal structure on $Y$ is $X$.

One example of such a comparison is the following result which is a direct application of an inequality of G. Anderson (\cite[Theorem 4.2]{anderson-thesis}) comparing the hyperbolic and projective metrics on a projective structure.
\begin{theorem}[Anderson, {\cite[Theorem 2.10]{BBB}}]
If $\Sigma$ is a projective structure with Schwarzian parameterization $\left(X_\Sigma, \Phi_\Sigma\right)$ and Thurston parameterization $\left(Y_\Sigma, \lambda_\Sigma\right)$ then
$$L\left(\lambda_\Sigma\right) \le 4\pi |\chi(\Sigma)|\left\|\Phi_\Sigma\right\|_\infty.$$
\end{theorem}
In \cite{BBB} bounds on $L\left(\lambda_\Sigma\right)$ in terms of the $L^2$-norm of $\Phi_\Sigma$ are also obtained.

The goal of this note is to bound  the $L^2$-norm of the quadratic differential in terms of the length of the lamination. We prove the following:
\begin{theorem}\label{main}
 Let $\Sigma$ be a projective structure that is the quotient of a domain in $\chat$ with Schwarzian parameterization $\left(X_\Sigma, \Phi_\Sigma\right)$ and Thurston parameterization $\left(Y_\Sigma, \lambda_\Sigma\right)$. Then
$$\left\|\Phi_\Sigma\right\|_2 \leq \left(1+\left\|\Phi_\Sigma\right\|_\infty\right)\sqrt{L\left(\lambda_\Sigma\right)} .$$
\end{theorem}
The restriction to projective structures that are quotients of domains in $\chat$ is not necessary but it will simplify the proof. On the other hand, if $\Sigma$ is the quotient of a disk in $\chat$ then by the Nehari bound on the Schwarzian (\cite{nehari}) $\left\|\Phi_\Sigma\right\|_\infty \le 3/2$ and we have the immediate corollary:
\begin{corollary}\label{maincor}
Let $\Sigma$ be a projective structure that is the quotient of a disk in $\chat$.  Then 
$$\left\|\Phi_\Sigma\right\|_2 \leq \frac{5}{2}\sqrt{L\left(\lambda_\Sigma\right)} .$$
\end{corollary}

Thurston's parameterization of projective structures involves associating the projective structure with a locally convex pleated surface in $\Hs$, pleated along the lamination $\lambda$. A parallel construction of Epstein gives an immersed surface whose geometry is described by the Schwarzian $\Phi$. The proof of Theorem \ref{main} involves the {\em $W$-volume} of the region between the two surfaces. This is a concept introduced by Krasnov and Schlenker \cite{KS08}. While it is generally applied to convex, compact hyperbolic 3-manifolds there is a natural relative version for pairs of conformal metrics on a projective structure and this is what we will use.

\section{Epstein Surfaces}
In \cite{epstein-envelopes}, Epstein described how to obtain a surface in $\Hs$ from a conformal metric on a domain in $\partial \Hs = \hat\CC$. Let $\hat g$ be a conformal metric on a  domain $\Omega$ of $\hat\CC$. For each $x \in \Hs$ we define the $\rho_x$ to be the {\em visual metric} on $\hat\CC$ from $x$. Thus $\rho_x$ is the image of the round metric on the sphere under any mobius transformation sending $0$ to $x$ in the Poincare ball model of $\Hs$. Given $z \in \Omega$ we define
$$\mathfrak H(\hat g,z) = \left\{ x\in \Hs\  | \ \hat g(z) = \rho_x(z) \right\}$$
Then $\mathfrak H(\hat g, z)$ is a horosphere based at $z$. Then the  {\em Epstein map} is a map
$$\Ep_{\hat g}\colon \Omega \to \Hs$$
such that for each $z\in \Omega$ the image surface is tangent to the horosphere $\mathfrak H(\hat g, z)$. 
%

We now collect some fundamental properties of Epstein surfaces and their shape operators that we will need.

\begin{prop}[{Epstein \cite{epstein-envelopes}}]\label{epstein}
Let $\hat g$ be a smooth conformal metric on a domain $\Omega \subset \chat$ and let $\hat g_t = e^{2t} \hat g$. 
\begin{enumerate}
\item For each $z \in \Omega$ the map $\Ep_{\hat g_t}$ is an immersion in a neighborhood of $z$ for all but at most two values of $t$. 
\item The map $\Ep_{\hat g_t}$ is the time $t$ of the normal flow of $\Ep_{\hat g}$. 
\end{enumerate}
When $\Ep_{\hat g_t}$ is an immersion let $g_t$ is the pullback via $\Ep^*_{\hat g_t}$ of the hyperbolic metric on $\Hs$ and  $B_t$  the shape operator for the immersed surface. Then
\begin{enumerate}
\setcounter{enumi}{2}
\item $\hat g_t = (\Id + B_t)^* g_t;$

\item $B_t = (\cosh (t-s) \cdot \Id + \sinh (t-s) \cdot B_s)^{-1}(\sinh (t-s)\cdot \Id + \cosh (t-s) \cdot B_s);$
\item $-1$ is never an eigenvalue of $B_t$.
\end{enumerate}
If we let $\hat B_t = (\Id + B_t)^{-1}(\Id - B_t)$ then
\begin{enumerate}
\setcounter{enumi}{5}
\item $g_t = \frac14(\Id + \hat B_t)^*\hat g_t$;

\item $\hat B_t = e^{-2(t-s)} \hat B_s$.
\end{enumerate}
The eigenvalues of $B_t$ are non-negative if and only if the eigenvalues of $\hat B_t$ are in the interval $(-1,1]$.
\end{prop}

 Property (2) is \cite[Theorem 2.1]{epstein-envelopes} while  (3) and  (4) are  direct consequences of  \cite[Theorem 3.1]{epstein-envelopes}. The operator $\hat B_t$ was introduced by Krasnov-Schlenker (\cite{KS08}). The remaining properties follow immediately from the others.

Note the formula in (4) only makes sense at a point where $B_s$ is defined and when $(\cosh (t-s)\cdot \Id + \sinh (t-s)\cdot B_s)^{-1}$ is invertible. If $B_s$ is defined then this holds as long as none of the eigenvalues of $B_s$ are $\coth (t-s)$. As $B_s$ has two eigenvalues (counted with multiplicity) this can happen at most twice and this explains (1) that at any point in $\Omega$ the map $\Ep_{\hat g_t}$ is an immersion for all but at most two values of $t$. By (5), whenever $B_t$ is defined so is $\hat B_t$. However, from (7) it is clear that $\hat B_t$ extends continuously and smoothly to all $t$. Note that all the information in $\hat B_t$ is contained in $B_t$ but this fact, along with the simplicity of the formula, are the advantage of considering $\hat B_t$.

\subsubsection*{The hyperbolic metric}
If $\hat g_h$ is the hyperbolic metric for $\Omega$ then the shape operator $B$ for the Epstein surface has a very explicit description in terms of the {\em Schwarzian derivative} of uniformizing map. In particular,  if $\chat \smallsetminus \Omega$ contains at least 3 points there is a  universal covering map
$$f\colon \Delta \to \Omega$$
from the unit disk $\Delta$ that is conformal. This map is unique up to pre-composition with M\"obius transformations fixing the disk and its Schwarzian derivative determines a holomorphic quadratic differential $\Phi_\Omega = \phi_\Omega dz^2$ on $\Omega$. Then $\hat g_h = \rho_h \hat g_{\rm euc}$ where $\rho_h$ is a smooth function and $\hat g_{\rm euc}$ is the Euclidean metric. The pointwise norm of the Schwarzian is given by the function $\|\Phi_\Omega(z)\| = \frac{|\phi_\Omega(z)|}{\rho_h(z)}$. We let $\|\Phi_\Omega\|_\infty$ be the $L^\infty$-norm of this function.

\begin{prop}[{Epstein, \cite[Proposition 7.4]{epstein-envelopes}}]\label{hyperbolic}
If $\hat g_h$ is the hyperbolic metric for $\Omega$ then the eigenvalues of the shape operator $B$ at $z\in\Omega$ are $-\frac{\|\Phi_\Omega(z)\|}{\|\Phi_\Omega(z)\|\pm 1}$. Therefore the eigenvalues of $\hat B_t$ are $e^{-2t}(1\pm 2\|\Phi_\Omega(z)\|)$ and the shape operator $B_t$ has non-negative eigenvalues if $e^{2t} \ge 1+ 2\|\Phi_\Omega\|_\infty$.
\end{prop}

\subsubsection*{Convex Epstein surfaces}
For certain surfaces in $\Hs$ we can also reverse this process (whenever the {\em hyperbolic Gauss map} is injective). The setting that will be important here is when the surface bounds a convex region in $\Hs$. Let $C\subset \Hs$ be a closed convex region and $\bar C$ the closure of $C$ in $\Hs\cup \chat$. Then $\Lambda = \bar C \cap \chat$ is a closed subset of $\chat$. Then $\Omega = \chat \smallsetminus \Lambda$ is open and for each $z\in \Omega$ there will be a unique horosphere that meets $C$ at a single point $x$ and we define $\hat g_C = \rho_x$ at $z$.  If $\del C$ is smooth then $\hat g_C$ will be smooth and  we have $\Ep_{\hat g_C}(z) = x$. 

If the boundary of $C$ is not smooth then the metric $\hat g$ is still defined and it will be continuous but not necessarily   differentiable. We can take $\Ep_{\hat g_C}(z) = x$ as the definition of the Epstein map. The map will be continuous and surjective but not necessarily smooth or injective. As the map is not smooth we cannot pull back the hyperbolic metric on $\Hs$ to $\Omega$. We will only be interested in a particular special case where $\del C$ is not smooth. We describe this next.

Let ${\rm CH}(\Lambda)$ be the smallest closed, convex subset of $\Hs$ such that $\overline{{\rm CH}(\Lambda)} \cap \chat = \Lambda$.
%
The boundary of ${\rm CH}(\Lambda)$ is a {\em convex pleated surface}.
In particular, the intrinsic path metric on $\del {\rm CH}(\Lambda)$ is a complete hyperbolic metric and in a complement of a geodesic lamination $\lambda$ the surface is totally geodesic. The surface is bent along $\lambda$ and the amount of bending gives $\lambda$ a transverse measure making it a measured geodesic lamination.

There is another intrinsic description of $\hat g_{{\rm CH}(\Lambda)}$.
First,  let $\hat g_D$ be the hyperbolic metric on a round disk $D \subset \chat$. Then for each $z\in \Omega$ define $\hat g_\Omega$ at $z$ to be the infimum of $\hat g_D$ at $z$ over all $D\subset \Omega$. This is the {\em projective metric} for $\Omega$ and $\hat g_\Omega = \hat g_{{\rm CH}(\Lambda)}$.

\section{Complex projective structures}
A {\em complex projective structure} $\Sigma$ on an surface $S$ is an atlas to $\chat$ where the transition maps are restrictions of elements of $\psl$. One example is when $\Omega$ is a domain in $\chat$ and $\Gamma \subset \psl$ is a deck acton on $\Omega$. Then the quotient $\Sigma = \Omega/\Gamma$ has a natural complex projective structure.  As the Epstein map is natural, if $\hat g$ is a $\Gamma$-invariant conformal metric on $\Omega$ the $\Ep_{\hat g}$ will also be $\Gamma$-invariant and the $\Ep_{\hat g}$-pullback $g$ of the hyperbolic metric on $\Hs$ will also be $\Gamma$-invariant. We will abuse notation and refer to the quotient metrics on $\Sigma$ also as $\hat g$ and $g$. Note that $g$ may be singular but if the image is a convex surface $g$ will be a smooth metric.

\begin{prop}[{{\cite[Lemma 3.1]{BBB}}}]\label{mean_integral}
Let $\hat g$ be a smooth conformal metric on a projective structure $\Sigma = \Omega/\Gamma$ and assume that the eigenvalues of the shape operator $B$ for Epstein surface are non-negative. If $g$ is the Epstein metric then
$$\int_S H dA = \frac12\area(\hat g) -\area(g) - \pi\chi(S)$$
where $H = \frac 12 \tr B$ is the mean curvature.
\end{prop}

Now let $\tilde C \subset \Hs$ be a $\Gamma$-invariant closed, convex subset  such that intersection of the closure of $\tilde C$ with $\chat$ is $\chat \smallsetminus  \Omega$. While the boundary of $C = \tilde C/\Gamma$ may not be smooth we can still take the area of the induced path metric $g_C$ and of the projective metric $\hat g_C$ on $\Sigma$. In particular, we can define
$$\cB(C) = \frac 12\area(\hat g_C) -\area(g_C) -\pi\chi(\Sigma).$$

Let $\hat g_h$ be the hyperbolic metric for the conformal structure on $\Sigma$ and let $\hat g_\Sigma$ be the projective metric. By Proposition \ref{hyperbolic} the shape operators for the Epstein metric of $(\hat g_h)_t = e^{2t}\hat g_h$ have non-negative eigenvalues when $e^{2t}\ge 1 +2\|\Phi_\Omega\|_\infty$. For the projective metric, the rescaled metrics $(\hat g_\Sigma)_t = e^{2t}\hat g_\Sigma$ are the conformal metric of a convex subset whenever $t\ge 0$.  A key to our estimate is that in both case we can explicitly calculate the area of the metrics at infinity and the corresponding Epstein metric.

For the hyperbolic metric we have:
\begin{lemma}\label{intmean}
Let $\hat g_h$ be the hyperbolic metric  on $\Sigma$ and $(g_h)_t$ be the intrinsic metric on $\Ep_{e^{2t}\hat g_h}$. Then  for $t >\frac12\log \left(1+2\left\|\Phi_\Sigma\right\|_\infty\right)$ we have
$$\area((g_h)_t) = -2\pi\chi(\Sigma) \cosh^2(t)-e^{-2t}\left\|\Phi_\Sigma\right\|^2_2 .$$

\end{lemma}

{\bf Proof:}
If $d{\hat A}_t$ is the area form for $(\hat g_h)_t$ and $dA_t$ is the area form for $g_t$ then by (4) and (5) of Proposition \ref{epstein} we have
$$d A_t = \frac14\left|\det\left(\Id + \hat B_t\right)\right| d\hat A_t= \frac14\left|\det\left(\Id + e^{-2t} \hat B_0\right)\right| e^{2t}d\hat A_0.$$
The determinant of $\Id + \hat B_t$ is non-negative when the eigenvalues of $\hat B_t$ are $\ge -1$ and by Proposition \ref{hyperbolic} this is exactly when $e^{2t} \ge 1 + 2\left\|\Phi_\Sigma\right\|_\infty$. Therefore when  $t >\frac12\log \left(1+2\left\|\Phi_\Sigma\right\|_\infty\right)$ we can drop the absolute value signs and integrate to get
$$\area((g_h)_t) =\int dA_t = \frac{1}{4}\int \det\left(\Id+e^{-2t}\hat B_0\right) e^{2t} d\hat A_0.$$
By  Proposition \ref{hyperbolic}, $\hat B$ has eigenvalues $1\pm 2\|\phi_\Sigma(z)\|$.  Thus
$$\det\left(\Id+e^{-2t}\hat B_0\right) = \left(1+e^{-2t}(1+ 2\left\|\Phi_\Sigma\right\|)\right)\left(1+e^{-2t}(1- 2\left\|\Phi_\Sigma\right\|)\right)=( 1+e^{-2t})^2 - 4\left\|\Phi_\Sigma\right\|^2e^{-4t}.$$
As $\area(\hat g_h) = -2\pi\chi(\Sigma)$ we have
$$\area((g_h)_t) =  \frac{1}{4}\left(e^{t}+e^{-t})^2 (-2\pi\chi(\Sigma)\right) -e^{-2t}\left\|\Phi_\Sigma\right\|^2_2 = -2\pi\chi(\Sigma)\cosh^2(t) -e^{-2t}\left\|\Phi_\Sigma\right\|^2_2.$$
\eproof

For the projective metric $\hat g_\Sigma$ we have the following;

\begin{lemma}\label{area_pr}
Let $\Sigma$ be a projective structure with Thurston parameterization $\left(Y_\Sigma, \lambda_\Sigma\right)$. Then
$$\area(\hat g_\Sigma) = -2\pi\chi(\Sigma) + L\left(\lambda_\Sigma\right)$$
and
$$\area((g_\Sigma)_t) = -2\pi\cosh^2(t)\chi(\Sigma) + L\left(\lambda_\Sigma\right)\sinh(t)\cosh(t).$$
\end{lemma}

{\bf Proof:} In the Thurston parameterization of projective structures $\Sigma$ is given by a pair $\left(Y_\Sigma, \lambda_\Sigma\right)$ (see \cite{Kamishima:Tan} for a complete description of Thurston's parametrization). We can assume that the support of $\lambda_\Sigma$ is a single curve $\gamma$ with weight $\theta$ so that $L(\lambda_\Sigma) = \theta\ell_{Y_\Sigma}(\gamma)$. The Epstein surface at time $t$ is the union of a Euclidean cylinder of circumference $\ell_{Y_\Sigma}(\gamma)\cosh(t)$ and width $\theta\sinh(t)$ with the complement a surface of constant negative curvature equal to $\tanh^2(t)-1$. The area of the cylinder is $\theta\ell_{Y_\Sigma}(\gamma) \sinh(t)\cosh(t) = L(\lambda_\Sigma)\sinh(t)\cosh(t)$. By the Gauss-Bonnet formula the area of the remainder of the surface is $-2\pi\cosh^2(t)\chi(\Sigma)$. This gives the formula for $\area((g_\Sigma)_t)$.

A similar computation gives the $\area(\hat g_\Sigma)$. In this metric the annulus has circumference $\ell_{Y_\Sigma}(\gamma)$ and width $\theta$ and the complementary surface is hyperbolic. The reasoning of the previous paragraph implies the formula for $\area(\hat g_\Sigma)$.

Laminations supported on single curve are dense in the space of all measure laminations and the lengths and areas will vary continuously. This gives the general statement. \eproof
%
\section{$W$-volume}
We first describe $W$-volume for convex cocompact hyperbolic manifolds and then its generalization for projective structures. Let $M = \Hs/\Gamma$ be a convex co-compact hyperbolic three-manifold with domain of discontinuity $\Omega$, limit set $\Lambda = \chat - \Omega$ and  conformal boundary  $\partial_c M = \Omega/\Gamma$. As $\Gamma$ acts by M\"obius transformations, we have a natural projective structure $\Sigma$ on $\partial_c M$.

Let $\tilde N$ be a closed, $\Gamma$-invariant, convex subset of $\Hs$ and assume that $N = \tilde N/\Gamma$ is compact. This is equivalent to the closure of $\tilde N$ in $\bar{\mathbb H}^3$  intersects $\chat$ in $\Lambda$. The $W$-volume  of $N$ is 
$$W(N) = \vol(N) -\frac12 \cB(N).$$
%
%

In \cite{KS08}, Krasnov and Schlenker define the $W$-volume by subtracting one half the integral of the mean curvature of the boundary. By Proposition \ref{mean_integral}, when $\del N$ is smooth this is an equivalent definition. If we gave the space of convex submanifolds the Gromov-Hausdorff topology $W$-volume will be continuous and convex submanifolds with smooth boundary will be dense.

The $W$-volume has a number of nice analytic properties. 
Let $N_t$ be the $t$-neighborhood of $N$.  Then a simple calculation (see   \cite[Lemma 4.2]{KS08}) shows that the $W$-volume satisfies the scaling property:
\begin{equation}\label{scaling}
W(N_t) = W(N) + t\pi|\chi(\partial M)|.
\end{equation}

Let $\cM_C(\del_c M)$ be conformal metrics $\hat g$ on $\del_c M$ such that
  there exists a $t$ for which the Epstein surface for $e^{2t} \hat g$  is the metric at infinity for a convex submanifold $N_t \subset M$.
 Using the scaling formula we can then define $W$-volume as a function on $\cM_C(\del_c M)$ by setting
$$W(\hat g) =  W(N_t) -t\pi |\chi(\del M)|.$$
By \eqref{scaling} the definition of $W(\hat g)$ doesn't depend on the choice of $t$.


\begin{theorem}[{Krasnov-Schlenker, \cite[Section 7]{KS08}}]\label{max_cc}
If $\hat g_{M}$ is the hyperbolic metric on $\del_c M$ and $\hat g \in \cM_C(\del_c M)$ has the same area as $\hat g_{M}$ then $W(\hat g_{M}) \ge W(\hat g)$.
\end{theorem}

%
%
%
%
%
%
%

\newcommand{\GG}{{\mathcal G}}
\newcommand{\PP}{{\mathcal P}}



\subsection{$W$-volume for projective structures}
We now define a relative version of $W$-volume and show that it naturally extends to a pair of conformal metrics on a projective structure. 

For a convex co-compact hyperbolic 3-manifold $M$ with metrics $\hat g_0, \hat g_1 \in \cM_C(\del_c M)$ we define
$$W(\hat g_0, \hat g_1) = W(\hat g_1) - W(\hat g_0)$$
and note that if the $\hat g_i$ are the metrics at infinity for convex submanifolds $N_i \subset M$ with $N_0 \subset N_1$ then
$$W(\hat g_0, \hat g_1) = \vol(N_1\smallsetminus N_0) - \frac12\left(\cB(N_1) - \cB(N_0)\right).$$

%
%

With this as motivation we now define a relative $W$-volume for pairs of conformal metrics on a projective structure. The construction can be made to work for any projective structure. However, for simplicity we'll restrict to project structures that are the quotient of domains in $\chat$. 

Here is our setup. Let $\Omega$ be a domain in $\chat$ and $\Gamma \subset \psl$ a deck action on $\Omega$ such that the projective structure $\Sigma = \Omega/\Gamma$ is compact. Let $\Lambda = \chat \smallsetminus \Omega$. As $\Gamma$ is a deck action on $\Omega$ it must be a discrete and torsion free subgroup of $\psl$, so the quotient $M = \Hs/\Gamma$ is a hyperbolic 3-manifold. If $\tilde N \subset \Hs$ is a closed $\Gamma$-invariant, convex set whose closure in $\chat$ is $\Lambda$ then the quotient $N = \tilde N/\Gamma$ will be a convex submanifold of $M$. However, $N$ will not necessarily be compact and the volume of $N$ may be infinite.

Now take nested convex, submanifolds $N_0 \subset N_1$ in $M$ such that the closure of the universal cover of the $N_i$ in $\chat$ is $\Lambda$. Then $N_1\smallsetminus N_0$ will have compact closure in $M$. If $\hat g_0$ and $\hat g_1$ are the metrics at infinity for $N_0$ and $N_1$ define
$$W\left(\hat g_0, \hat g_1\right) = \vol(N_1\smallsetminus N_0) - \frac12\left(\cB(N_1) - \cB(N_0)\right)$$
as above. Recalling that the metric at infinity for the $t$-neighborhood of $N_1$ is $e^{2t} \hat g_1$ as for the $W$-volume for convex co-compact hyperbolic 3-manifolds we have the scaling property
$$W\left(\hat g_0, e^{2t}\hat g_1\right) = W\left(\hat g_0, \hat g_1\right) + t\pi|\chi(\Sigma)|.$$
As before this formula only makes sense initially when $t \ge 0$. For $t<0$ we can take it as a definition of $W\left(\hat g_0, e^{2t}\hat g_1\right)$. The scaling property gives that this is well defined. If the $s$-neighborhood of $N_0$ is contained in $N_1$ then the scaling property gives
$$W\left(e^{2s}\hat g_0, \hat g_1\right) = W\left(\hat g_0, \hat g_1\right) -s\pi|\chi(\Sigma)|.$$
Again we take this as the definition of $W\left(\hat g_0, e^{2s}\hat g_1\right)$ for arbitrary $s$. Finally if $\hat g_0$ and $\hat g_1$ are arbitrary conformal metrics on $\Sigma$ that can be rescaled to be the metrics at infinity for convex submanifolds of $M$ we can choose the scaling factors $t$ and $s$ such that $e^{2t}\hat g_1 \ge e^{2s}\hat g_0$ and define
$$W(\hat g_0, \hat g_1) = W\left(e^{2t}\hat g_0, e^{2s}\hat g_1\right) -(t-s)\pi|\chi(\Sigma)|.$$
It then follows that
$$W(\hat g_0, \hat g_1) = -W(\hat g_1, \hat g_0).$$
%
%
%

As in the case of a convex co-compact hyperbolic 3-manifold we let $\cM_C(\Sigma)$ be the set of conformal metrics $\hat g$ on $\Sigma$ such that there exists a $t$ with $e^{2t}\hat g$ the metric at infinity for some convex $N \subset M$ as above. We can define $W\left(\hat g_0, \hat g_1\right)$ for any pair of metric $\hat g_0, \hat g_1$ in $\cM_C(\Sigma)$. The relative version of the Theorem \ref{max_cc} becomes:
\begin{theorem}\label{max_pr}
If $\hat g_h$ is the hypebrolic metric for the conformal structure on $\Sigma = \Omega/\Gamma$ and $\hat g \in \cM_C(\Sigma)$ has the same area as $\hat g_h$ then
$$W(\hat g_h, \hat g) \le 0.$$
\end{theorem}

\section{Calculating $W(\hat g_h,\hat g_\Sigma)$}
The key to our estimate is that for both the hyperbolic metric and the projective metric on $\Sigma$ we can explicitly calculate the area of the surface at infinity and the Epstein surface and furthermore we can make this calculation after rescaling.

\begin{lemma}\label{upper}
$$W(\hat g_h, \hat g_\Sigma) \le \frac14 L\left(\lambda_\Sigma\right)$$
\end{lemma}

{\bf Proof:}
By Lemma \ref{area_pr}
$$\area\left(\hat g_\Sigma\right) = 2\pi|\chi(\Sigma)|+L(\lambda_\Sigma).$$
Applying the maximality property of Theorem \ref{max_pr} we have
$$W\left(\hat g_h,\frac{2\pi|\chi(\Sigma)|}{2\pi|\chi(\Sigma)|+L(\lambda_\Sigma)}\hat g_\Sigma\right) \leq 0.$$
Combing this with the scaling property we get
$$W(\hat g_h,\hat g_\Sigma) +\frac{1}{2}\log\left(\frac{2\pi|\chi(\Sigma)|}{2\pi|\chi(\Sigma)|+L(\lambda_\Sigma)}\right)\pi|\chi(\Sigma)| \leq 0.$$
As  $\log(1+x) \le x$ we have
$$W(\hat g_h,\hat g_\Sigma) \leq \frac{1}{2}\log\left(1+\frac{L(\lambda_\Sigma)}{2\pi|\chi(\Sigma)|}\right)\pi|\chi(\Sigma)| \leq \frac{1}{4}L(\lambda_\Sigma).$$
\eproof

Next we get a lower bound. The key here is that volumes are non-negative.

\begin{lemma}\label{lower}
$$W(\hat g_h, \hat g_\Sigma) \ge \frac{e^{-2T}\left\|\Phi_\Sigma\right\|_2^2}2 - \frac{L\left(\lambda_\Sigma\right)\cosh(2T)}4$$
where $e^{2T} = 1+2\left\|\Phi_\Sigma\right\|_\infty$.
\end{lemma}

{\bf Proof:} 
By Proposition \ref{hyperbolic} the eigenvalues for the shape operator for the Epstein surface of $(\hat g_h)_T$ are non-negative so $(\hat g_h)_T $ is the metric at infinity for a convex submanifold $(N_h)_T$.  For all $t\ge 0$, $(\hat g_\Sigma)_t$ is the metric at infinity for a convex submanifold $(N_\Sigma)_t$. Since $\hat g_h \le \hat g_\Sigma$ we have $(N_h)_T \subset (N_\Sigma)_T$ and therefore
\begin{eqnarray*}
W(\hat g_h, \hat g_\Sigma) & = & W\left((\hat g_h)_T, (\hat g_\Sigma)_T\right) \\
&= & \vol\left((N_\Sigma)_T\smallsetminus (N_h)_T\right) -\frac12\left(\cB\left((N_\Sigma)_T\right) - \cB\left((N_h)_T\right)\right)\\
&\geq& -\frac12\left(\cB\left((N_\Sigma)_T\right) - \cB\left((N_h)_T\right)\right)\\
&\ge & \frac12\left(\frac{L\left(\lambda_\Sigma\right)e^{2T}}2 + e^{-2T}\left\|\Phi_\Sigma\right\|_2^2 - L(\lambda_\Sigma)\sinh(T)\cosh(T)\right)\\
& \ge &\frac12 \left(e^{-2T}\left\|\Phi_\Sigma\right\|_2^2 -\frac{L(\lambda_\Sigma)}2 \cosh(2T)\right).
\end{eqnarray*}
Here we are applying Lemmas \ref{intmean} and \ref{area_pr} to calculate $\cB\left((N_h)_T\right)$ and $\cB\left((N_\Sigma)_T\right)$.
\eproof

We now prove our main result. \newline
\vspace{.001in}

\noindent
{\bf Theorem \ref{main}}
{\em Let $\Sigma$ be a projective structure that is the quotient of a domain in $\chat$ with Schwarzian parameterization $\left(X_\Sigma, \Phi_\Sigma\right)$ and Thurston parameterization $\left(Y_\Sigma, \lambda_\Sigma\right)$. Then
$$\left\|\Phi_\Sigma\right\|_2 \leq \left(1+\left\|\Phi_\Sigma\right\|_\infty\right)\sqrt{L\left(\lambda_\Sigma\right)} .$$

}
{\bf Proof:} Combining Lemmas \ref{upper} and \ref{lower} we get
$$\frac{e^{-2T}\left\|\Phi_\Sigma\right\|_2^2}2 - \frac{L\left(\lambda_\Sigma\right)\cosh(2T)}4\leq \frac{1}{4}L(\lambda_\Sigma)$$
giving
$$\left\|\Phi_\Sigma\right\|_2^2 \leq \frac{1}{2}L(\lambda_\Sigma)e^{2T}(1+\cosh(2T)) =   \frac{1}{4}L(\lambda_\Sigma)\left(e^{2T}+1\right)^2.$$
Applying $e^{2T} = 1+2\left\|\Phi_\Sigma\right\|_\infty$ we obtain 
$$\left\|\Phi_\Sigma\right\|_2 \leq   (1+\left\|\Phi_\Sigma\right\|_\infty)\sqrt{L(\lambda_\Sigma)}.$$
\eproof
\bibliography{bib,math}
\bibliographystyle{math}

\end{document}